\documentclass[titlepage,11pt]{article}
\usepackage{latexsym}
\usepackage{amsfonts}
\usepackage{amsmath}
\usepackage{mathtools}
\usepackage{amssymb}
\usepackage{tikz}
\usepackage{comment}
\newcommand\blackslug{\hbox{\hskip 1pt \vrule width 4pt height 8pt depth 1.5pt
		\hskip 1pt}}
\newcommand\bbox{\hfill \quad \blackslug \bigbreak}
\def\DD{\hbox{-}}
\def\CC{\hbox{-}\cdots\hbox{-}}
\def\LL{,\ldots,}

\newcommand{\dist}{\operatorname{dist}}

\newcommand{\mac}{\mathcal}

\title{Asymptotic structure. I. Coarse tree-width}
\author{
	Tung Nguyen\thanks{Part of this work was done while the first author was at Princeton University and was supported by a Porter Ogden Jacobus Fellowship, AFOSR grant
		FA9550-22-1-0234, and NSF grant  DMS-2154169. Currently supported by a Titchmarsh Research Fellowship and a Christ Church Research Centre Grant.}\\
	University of Oxford,\\
	Oxford, UK
	\and
	Alex Scott\thanks{Supported by EPSRC grant EP/X013642/1}\\
	University of Oxford, \\
	Oxford, UK
	\and
	Paul Seymour\thanks{Supported by AFOSR grant
		FA9550-22-1-0234, and by NSF grant DMS-2154169.}\\
	Princeton University,\\ Princeton, NJ 08544, USA}

\date{}

\newtheorem{thm}{}[section]

\newcommand{\Proof}{\noindent{\bf Proof.}\ \ }

\begin{document}
	\maketitle
	\begin{abstract}
		In this paper, we develop a coarse analogue of tree-width.
		We prove that a graph $G$ 
		admits a tree-decomposition in which each bag is contained in the union of a bounded number of balls of bounded radius, if and only if 
		$G$ admits a quasi-isometry to a graph with bounded tree-width. (The ``if'' half is easy, but the 
		``only if'' half is challenging.) This generalizes a recent result of Berger and Seymour, concerning tree-decompositions when 
		each bag has bounded radius. We also prove a similar result for line-width, which is an extension of path-width to infinite graphs.
		
	\end{abstract}

	\section{Introduction}
	
	This is the first in a sequence of papers on asymptotic structure in graphs, looking at properties of graphs that are preserved under 
	quasi-isometry.  The notion of quasi-isometry is well established in geometric group theory: following work of Gromov 
	(see~\cite{gromov1,gromov2}), the large-scale geometric structure of a
	group can be investigated by looking at properties of its 
	Cayley graphs that are invariant under quasi-isometry.  More recently, the emerging area of ``coarse graph theory'' has begun to 
	look at large-scale geometry of graphs in general, using the lens of quasi-isometry.  It has begun to emerge that many classical 
	notions from graph theory have counterparts in the coarse context: for example, a coarse analogue for graph minors is provided by 
	``fat minors'' (see~\cite{jems, chepoi}); and there are many questions about which graph-theoretic results have natural counterparts 
	in the coarse context~\cite{gp25}.  In this paper, we investigate a coarse analogue of tree-width that was introduced in~\cite{dragan}, and connect it via quasi-isometry to the 
	standard notion of tree-width.  Our results also hold for path-width; however, more can be said about path-width, and we will pursue 
	this further in~\cite{AS2, AS3}.
	
	We need to begin with some definitions.
	Graphs in this paper may be infinite.
	If $X$ is a vertex of a graph $G$, or  
	a subset of the vertex set of $G$, or a subgraph of $G$, and the same for $Y$, then $\dist_G(X,Y)$ denotes the 
	distance in $G$
	between $X,Y$, that is, the number of edges in the shortest path of $G$ with one end in $X$ and the other in $Y$. (If no path exists we set $\dist_G(X,Y) = \infty$.)
	
	Let $G,H$ be graphs, and let $\phi:V(G)\to V(H)$ be a map.
	Let $L,C\ge 0$; we say that $\phi$ is an {\em $(L,C)$-quasi-isometry}
	if:
	\begin{itemize}
		\item for all $u,v$ in $V(G)$, if $\dist_G(u,v)$ is finite then $\dist_H(\phi(u),\phi(v))\le L \dist_G(u,v)+C$;
		\item for all $u,v$ in $V(G)$, if $\dist_H(\phi(u),\phi(v))$ is finite then $\dist_G(u,v)\le L \dist_H(\phi(u),\phi(v))+C$;
		and
		\item for every $y\in V(H)$ there exists $v\in V(G)$ such that $\dist_H(\phi(v), y)\le  C$.
	\end{itemize}
	%(See Manning~\cite{manning}.)
	If $X\subseteq V(G)$, let us say the {\em diameter of $X$ in $G$} is the maximum of $\dist_G(u,v)$ over all $u,v\in X$.
	A {\em tree-decomposition} of a graph $G$ is a pair $(T,(B_t:t\in V(T)))$, where $T$ is a tree (possibly infinite), and $B_t$
	is a subset of $V(G)$ for each $t\in V(T)$ (called a {\em bag}), such that:
	\begin{itemize}
		\item $G$ is the union of the graphs $G[B_t]\;(t\in V(T))$; and 
		\item for all $t_1,t_2,t_3\in V(T)$, if $t_2$ lies on the path of $T$ between $t_1,t_3$, then
		$B_{t_1}\cap B_{t_3}\subseteq B_{t_2}$.
	\end{itemize}
	The {\em width} of a tree-decomposition $(T,(B_t:t\in V(T)))$ is the maximum of the numbers $|B_t|-1$ for $t\in V(T)$,
	or $\infty$ if there is no finite maximum;
	and the {\em tree-width} of $G$ is the minimum width of a tree-decomposition of $G$.
	If $T$ is a path, we call $(T,(B_t:t\in V(T)))$ a {\em path-decomposition}, and the {\em path-width} of $G$ is defined analogously.

	Our first result is an extension of a result of Berger and Seymour~\cite{shortcuts} (which can also be derived from a combination 
	of results of Chepoi et al.~\cite{oldchepoi}). They proved:
	
	\begin{thm}\label{treeeqnce}
		For all $r$, if $G$ is connected and admits a tree-decomposition $(T,(B_t:t\in V(T)))$ such that for each $t\in V(T)$, $B_t$ has diameter at most $r$ in $G$,
		then $G$ admits a $(1,6r+1)$-quasi-isometry to a tree.
	\end{thm}
	This has a sort of converse, also proved in~\cite{shortcuts}: if $G$ is connected and $(L,C)$-quasi-isometric to a tree then it admits a tree-decomposition 
	$(T,(B_t:t\in V(T)))$ such that $B_t$ has diameter at most $L(L+C+1)+C$ in $G$, for each $t\in V(T)$.
	
	We will extend \ref{treeeqnce} from trees to graphs of bounded tree-width, as follows (although saying that this extends \ref{treeeqnce} is something of a stretch, because we
	do not know whether \ref{quasitw} holds with $L=1$):
	\begin{thm} \label{quasitw}
		For all $k,r$, there exist $L,C\ge 1$ such that if $G$ admits a tree-decomposition $(T,(B_t:t\in V(T)))$ such that for each $t\in V(T)$,
		$B_t$ is the union of at most $k$ sets each with diameter at most $r$ in $G$, then
		$G$ admits an $(L,C)$-quasi-isometry to a graph with tree-width at most $k$.
	\end{thm}
	A similar result (with weaker constants) was obtained independently by Hickingbotham~\cite{hick}, by applying a result of
	Dvo\v{r}\'{a}k and Norin~\cite{zdenek}.
	
	Our proof obtains a quasi-isometry to a graph with a tree-decomposition indexed by a subdivision of the same tree $T$ that 
	indexed the tree-decomposition of $G$; and so if $T$ is a path, we find a quasi-isometry to a graph with bounded path-width.
	Consequently:
	\begin{thm} \label{quasipw}
		For all $k,r$, there exist $L,C\ge 1$ such that if $G$ admits a path-decomposition $(T,(B_t:t\in V(T)))$ such that for each $t\in V(T)$,
		$B_t$ is the union of at most $k$ sets each with diameter at most $r$ in $G$, then
		$G$ admits an $(L,C)$-quasi-isometry to a graph with path-width at most $k$.
	\end{thm}
	
	In fact, for path-decompositions (and for ``line-decompositions'', discussed below), we can do much more: we can get an  {\em additive} 
	quasi-isometry to a graph with bounded path-width (or ``line-width''); that is, we can take $L=1$.  This follows from \ref{quasipw}, and the following result from~\cite{AS2}:
	
	\begin{thm}
		For all $L, C, k$ there exists $C'$
		such that if there is an $(L, C)$-quasi-isometry from a graph $G$ to
		a graph $H$ with line-width at most $k$, then there is a $(1, C')$-quasi-isometry from $G$ to a graph $H_0$ 
		obtained from $H$ by subdividing and contracting edges.    
	\end{thm}

	A path-decomposition is essentially a sequence of sets of vertices satisfying the ``betweenness'' condition. There is a more
	general notion (see~\cite{subtrees,diestel, thomas2, AS2,AS3}), where we replace the sequence by a family of subsets indexed by a linearly ordered set,
	giving what we call ``line-width''. Line-width and path-width are the same for finite graphs, but for infinite graphs they may be
	different. 
	We will show in the second half of the paper that \ref{quasipw} works with line-width in place of path-width. The proof uses the same basic approach 
	as that for \ref{quasitw}, but is significantly different in the details.

	How sharp is our bound on tree-width in \ref{treeeqnce}?
	In \ref{quasitw}, we start with a tree-decomposition in which each bag is the 
	union of $k$ bounded-radius balls, and we obtain a tree-decomposition in which each bag has size at most $k+1$: and one might hope that
	the final
	$k$ in the statement of \ref{quasitw} should be $k-1$.
	Obviously not for $k = 1$; but not when $k\ge 2$ either. To see this when $k=2$, let $G$ be a cycle, with vertices $v_1\CC v_n\DD v_1$
	in order. For $1\le i\le n-1$, let $B_{v_i} = \{v_i, v_{i+1}, v_n\}$, and let $T$ be the tree $G\setminus \{v_n\}$. 
	Then $(T, (B_t:t\in V(T)))$ is a tree-decomposition of $G$, and each of its bags is the union of two balls of bounded radius (one the singleton $\{v_n\}$ and the other consisting of two adjacent vertices).
	On the other hand, for all $(L,C)$, if $n$ is large enough then there is no $(L,C)$-quasi-isometry from $G$ to a graph 
	with tree-width at most $1$. A similar example works for each value of $k\ge 2$ (take a $k\times k$ grid and subdivide each of its edges many times).

	Again, \ref{quasitw} has a sort of converse, because if $G$ admits an $(L,C)$-quasi-isometry to a graph with tree-width at most $k$,
	then $G$ admits a tree-decomposition $(T,(B_t:t\in V(T)))$ such that for each $t\in V(T)$, $B_t$ is the union of at most
	$k+1$ sets each of bounded diameter --- we will prove this in the next section. But if we start with a graph $G$ that admits a quasi-isometry to a graph with tree-width at most $k$,
	and apply this converse, we obtain a tree-decomposition in which each bag is a union of $k+1$ sets of bounded diameter; and if we then apply
	\ref{quasitw}, we obtain a quasi-isometry to a graph with tree-width at most $k+1$. Somewhere we went from tree-width $k$ to tree-width $k+1$, and this is unsatisfying, at least on aesthetic grounds.
	
	A way to get rid of it is to make a small tweak in the definition of tree-decomposition; 
	say a {\em pseudo-tree-decomposition} $(T,(B_t:t\in V(T)))$ is the same as a tree-decomposition, except we relax the condition that every edge has both ends in some bag. 
	Instead, we insist that for every edge $uv$, either some bag contains both $u,v$, or there is an edge $st$ of $T$ such that 
	$B_s\setminus B_t = \{u\}$ and $B_t\setminus B_s = \{v\}$. Define {\em pseudo-tree-width} correspondingly (it differs from tree-width 
	by at most one). We will prove a version of \ref{quasitw} with ``tree-width at most $k$'' replaced by ``pseudo-tree-width at most $k$'',
	and a version of \ref{qtwconverse} with ``tree-width at most $k$'' replaced by ``pseudo-tree-width at most $k+1$'', and 
	the anomalous error of one is gone. More exactly, we will prove:
	
	\begin{thm}\label{pseudotw}
		For all $k,r$, there exist $L,C$ such that if $G$ admits a tree-decomposition $(T,(B_t:t\in V(T)))$ such that for each $t\in V(T)$,
		$B_t$ is the union of at most $k$ sets each with diameter at most $r$ in $G$, then
		$G$ admits an $(L,C)$-quasi-isometry to a graph with pseudo-tree-width at most $k$. 
		
		Conversely, for all $L,C\ge 1$, 
		if $G$ admits an $(L,C)$-quasi-isometry to a graph with pseudo-tree-width at most $k$,
		then $G$ admits a tree-decomposition $(T,(B_t:t\in V(T)))$ such that for each $t\in V(T)$, $B_t$ is the union of at most
		$k$ sets each of diameter at most $2L(L+C)+C$.
	\end{thm}
	
	Let us mention a related result.
	If a graph $G$ admits a quasi-isometry (with bounded constants) to a graph with tree-width at most $k$, then, trivially, there is a partition $\mac P$ of $V(G)$ such that:
	\begin{itemize}
		\item there is a graph $J$ with vertex set $\mac P$ such that $J$ has tree-width at most $k$;
		\item for each set $X\in \mac P$, every two vertices in $X$ have bounded distance in $G$;
		\item for every edge $XY$ of $J$, some vertex in $X$ has bounded distance in $G$ from some vertex in $Y$; and
		\item for all  distinct $X,Y\in \mac P$, if there is an edge of $G$ between $X,Y$ then the distance in $J$ between $X,Y\in V(J)$ is bounded.
	\end{itemize}
	A recent paper by Distel~\cite{distel} strengthens this: he shows that one can replace the third and fourth bullets above by:
	\begin{itemize}
		\item for all distinct $X,Y\in \mac P$, there is an edge of $G$ between $X$ and $Y$ if and only if $X,Y$ are adjacent in $J$.
	\end{itemize}
	As far as we know, this result of Distel is incomparable with our result; neither easily implies the other.

	%%%%%%%%%%%%%%%%%%%%%%%%%%%%%%%%%%%%%%%%%%%%%%%%%%%%%%%%%%%%%%%%%%%%%%%%%%%%%%%%%%%%%%%%%%%%%%%%%%%%%%%%%%%%%%%%%%%%%%%%%%

	\section{The proof of \ref{pseudotw}}
	
	Let us state the definition of pseudo-tree-width more formally.
	A {\em pseudo-tree-decomposition} of a graph $G$ is a pair $(T,(B_t:t\in V(T)))$, where $T$ is a tree, and $B_t$
	is a subset of $V(G)$ for each $t\in V(T)$ (called a {\em bag}), such that:
	\begin{itemize}
		\item $V(G)$ is the union of the sets $B_t\;(t\in V(T))$;
		\item for every edge $e=uv$ of $G$, either there exists $t\in V(T)$ with $u,v\in B_t$, or there is an edge $st\in E(T)$
		such that $B_s\setminus B_t = \{u\}$ and $B_t\setminus B_s = \{v\}$; and
		\item for all $t_1,t_2,t_3\in V(T)$, if $t_2$ lies on the path of $T$ between $t_1,t_3$, then
		$B_{t_1}\cap B_{t_3}\subseteq B_{t_2}$.
	\end{itemize}
	The {\em width} of a pseudo-tree-decomposition $(T,(B_t:t\in V(T)))$ is the maximum of the numbers $|B_t|$ for $t\in V(T)$,
	or $\infty$ if there is no finite maximum;
	and the {\em pseudo-tree-width} of $G$ is the minimum width of a pseudo-tree-decomposition of $G$. (Note that we are omitting the customary $-1$ in the definition of width.)
	If $T$ is a path, we call $(T,(B_t:t\in V(T)))$ a {\em pseudo-path-decomposition}, and the {\em pseudo-path-width} of $G$ is defined analogously.
	When $T$ is a finite path, we sometimes use the notation $(B_1\LL B_n)$ (with the usual meaning) in place of $(T,(B_t:t\in V(T)))$. 
	
	Before we prove the main part of \ref{pseudotw}, let us prove its (much easier) second part, the converse:
	\begin{thm}\label{qtwconverse}
		If $G$ admits an $(L,C)$-quasi-isometry to a graph with pseudo-tree-width at most $k$,
		then $G$ admits a tree-decomposition $(T,(D_t:t\in V(T)))$ such that for each $t\in V(T)$, $D_t$ is the union of at most
		$k$ sets each of diameter at most $2L(L+C)+C$.
	\end{thm}
	\Proof Let $H$ be a graph with pseudo-tree-width at most $k$, and let $(T,(B_t:t\in V(T)))$ be a pseudo-tree-decomposition of $H$ with
	width at most $k$. Let $\phi$ be an $(L,C)$-quasi-isometry from a graph $G$ to $H$.
	For each $h\in V(H)$, let $X_h$ be the set of vertices $i\in V(H)$ such that $\dist_H(h,i)\le L+C$.
	For each $t\in V(T)$, let $D_t$ be the set of all vertices $v\in V(G)$ such that $\phi(v)\in X_h$ for some $h\in B_t$.
	We claim that $(T,(D_t:t\in V(T)))$ is a tree-decomposition of $G$ satisfying the theorem. So we must check that:
	\begin{itemize}
		\item $\bigcup_{t\in V(T)}D_t = V(G)$;
		\item for every edge $uv$ of $G$ there exists $t\in V(T)$ with $\{u,v\}\in D_t$; 
		\item for all $t_1,t_2,t_3\in V(T)$, if $t_2$ lies on the path of $T$ between $t_1,t_3$, then $D_{t_1}\cap D_{t_3}\subseteq D_{t_2}$; and
		\item for each $t\in V(T)$, $D_t$ is the union of at most $k$ sets each of diameter (in $G$) at most $2L(L+C)+C$.
	\end{itemize}
	For the first statement, let $v\in V(G)$; then $\phi(v)\in V(H)$, and so $\phi(v)\in B_t$ for some $t\in V(T)$. In particular, since
	$\phi(v)\in X_{\phi(v)}$, it follows that $v\in D_t$. This proves the first statement.
	
	For the second statement, let $uv\in E(G)$, and choose $t\in V(T)$ with $\phi(v)\in B_t$. Since $\phi$ is an $(L,C)$-quasi-isometry,
	$\dist_H(\phi(u),\phi(v))\le L+C$,
	and so $\phi(u)\in X_{\phi(v)}$. It follows that $u,v\in D_t$. This proves the second statement.
	
	For the third statement, let $t_1,t_2,t_3\in V(T)$, such that $t_2$ lies on the path of $T$ between $t_1,t_3$, and 
	let $v\in D_{t_1}\cap D_{t_3}$. Hence for $i = 1,3$, there exists
	$h_i\in B_{t_i}$ with $\phi(v)\in X_{h_i}$;
	let $P_i$ be a path of $H$ between $\phi(v), h_i$ of length at most $L+C$.
	Since $P_1\cup P_3$ is a connected graph with vertices in $B_{t_1}$ and in $B_{t_3}$, it also has a vertex in $B_{t_2}$, say $h_2$.
	Thus $h_2$ belongs to one of $V(P_1), V(P_3)$, and so $\dist_H(h_2,\phi(v))\le L+C$; and hence $\phi(v)\in X_{h_2}$, and therefore
	$v\in D_{t_2}$. This proves the third statement.
	
	Finally, for the fourth statement, let $t\in V(T)$. For each $h\in B(t)$, let $F_h$ be the set of all $v\in V(G)$
	such that $\phi(v)\in X_h$. Thus $D_t$ is the union of the sets $F_h\;(h\in B_t)$, and there are $|B_t|\le k$ such sets.
	We claim that each $F_h$ has diameter at most $2L(L+C)+C$ in $G$. If $u,v\in F_h$, then each of $\phi(u), \phi(v)$ has distance at most $L+C$ from $h$,
	and so $\dist_H(\phi(u), \phi(v))\le 2(L+C)$. Since $\phi$ is an $(L,C)$-quasi-isometry, it follows that
	$\dist_H(u,v)\le 2L(L+C)+C$. This proves the fourth statement, and so proves \ref{qtwconverse}.~\bbox

	To prove \ref{pseudotw}, we need the following lemma:
	\begin{thm}\label{matching}
		Let $G$ be a graph, and let $A,B$ be disjoint subsets of $V(G)$ with union $V(G)$. Let $|A|,|B|\le k$, and suppose that there are
		at most $k$ edges between $A,B$. Then there is a pseudo-path-decomposition $(B_1\LL B_n)$ of $G$ with width at most $k$ and with $A\subseteq B_1$ 
		and $B\subseteq B_n$.
	\end{thm}
	\Proof
	We proceed by induction on $k+|A|+|B|$. If some vertex $a\in A$ has no neighbours in $B$, then from the inductive hypothesis,
	applied to $G\setminus \{a\}$, there is a pseudo-path-decomposition $(B_1\LL B_n)$ of $G\setminus \{a\}$ with width at most $k$ and with $A\setminus \{a\}\subseteq B_1$
	and $B\subseteq B_n$. But then $(A,B_1\LL B_n)$ satisfies the theorem. Thus we may assume that each vertex in $A$ has a neighbour in $B$, and vice versa.
	
	If every vertex in $A$ has exactly one neighbour in $B$ and vice versa, the result is true; so we assume that some vertex in $A$
	has at least two neighbours in $B$, and hence $|A|\le k-1$. Let $b\in B$ with a neighbour in $A$, and
	let $G'$ be obtained by deleting $b$. In $G'$, there are at most $k-1$
	edges between $A$ and $B\setminus \{b\}$, and these two sets both have size at most $k-1$.
	From the inductive hypothesis applied to $G'$, there is a 
	pseudo-path-decomposition $(C_1\LL C_n)$ of $G'$ with width at most $k-1$ and with $A\subseteq C_1$ and $B\setminus \{b\}\subseteq C_n$.
	Define $B_i = C_i\cup \{b\}$ for $1\le i\le n$; then $(B_1\LL B_n)$ is a pseudo-path-decomposition of $G$ satisfying 
	the theorem. This proves \ref{matching}.~\bbox

	To prove the first part of \ref{pseudotw}, it suffices to prove it when $G$ is connected (by working with each component of $G$ separately); and it suffices to prove it when $r=1$. To see the latter, let $G$ be a connected graph that admits a 
	tree-decomposition $(T,(B_t:t\in V(T)))$ such that for each $t\in V(T)$,
	$B_t$ is the union of at most $k$ sets each with diameter at most $r$ in $G$. For each $t\in V(T)$, and each pair $u,v$ of 
	nonadjacent vertices of $G[B_t]$ with 
	$\dist_G(u,v)\le r$, add an edge joining $u,v$, and let $G'$ be the resultant graph. Then 
	$(T,(B_t:t\in V(T)))$ is a tree-decomposition of $G'$, and for each $t\in V(T)$,
	$B_t$ is the union of at most $k$ cliques of $G'$. (Cliques are sets of vertices, not subgraphs.) 
	Moreover, the identity map is an $(r,0)$-quasi-isometry
	between $G,G'$; and so if $G'$ admits an $(L,C)$-quasi-isometry to a graph with pseudo-tree-width at most $k$, then $G$ admits an 
	$(rL,rC)$-quasi-isometry to the same graph. Consequently, for given $k$, if $L,C$ satisfy the theorem when $r=1$, then $rL,rC$ satisfy the theorem
	for general $r$. Hence it suffices to prove the result when $r=1$.
	
	We can weaken this hypothesis a little. If $B\subseteq V(G)$, $\alpha(B)$ denotes the size of the largest stable subset of 
	$G[B]$ (or $\infty$ if there is no such largest set). If we have a tree-decomposition in which each bag is the union of 
	at most $k$ cliques, then each bag $B$ satisfies $\alpha(B)\le k$, and this is enough for us. (Incidentally, the smallest $k$
	such that $G$ admits a tree-decomposition in which $\alpha(B)\le k$ for each bag $B$, is called the {\em tree-independence number} of $G$, and is of interest for algorithmic applications~\cite{treealpha1,treealpha2}.)
	We will prove the following:
	
	\begin{thm} \label{quasitw2}
		For all $k$, if $G$ is connected and admits a tree-decomposition $(T,(B_t:t\in V(T)))$ such that 
		$\alpha(B_t)\le k$ for each $t\in V(T)$, then
		$G$ admits a $(2k+2,2k-1)$-quasi-isometry to a graph with pseudo-tree-width at most $k$.
	\end{thm}
	\Proof Let $(T,(B_t:t\in V(T)))$ be a tree-decomposition of $G$ such that $\alpha(B_t)\le k$ for each $t\in V(T)$.
	Fix a root $r\in V(T)$ (arbitrarily). For each $t\in V(T)$, its {\em ancestors} are the vertices of the path
	of $T$ between $r,t$, and its {\em strict ancestors} are its ancestors different from $t$. If $s$ is an ancestor of $t$ then 
	$t$ is a {\em descendant} of $s$. If $s,t\in V(T)$, we write $s\le_T t$ if $s$ is an ancestor of $t$, and $s<_T t$ if $s$ is a strict ancestor of $t$.
	For $t\in V(T)$, its {\em height} is the length of the path of $T$ between $r,t$.
	
	We will recursively define a set of pairs, called ``cores''.
	Each core will be a pair $(t,C)$ where $t\in V(T)$ and $C$ is a subset of $B_t$ inducing a non-null connected subgraph, and we will call $t$
	its {\em birthday}. 
	Each core $(t,C)$ will have a {\em spread} 
	$S(t,C)$, which is the vertex set of a certain subtree of $T$ with root $t$, defined below.

	Here is the inductive definition.
	If there exists $t\in V(T)$ such that we have not yet defined the set of cores with birthday $t$,
	choose some such $t$ with minimum height.
	We say $v\in B_t$ is {\em disqualified} if there is a core $(s,C)$ such that $s<_T t$, and $t\in S(s,C)$,
	and either $v\in C$ or $v$ has a neighbour in $C$. Let $Z$ be the set of vertices in $B_t$ that are not disqualified. 
	For each $C\subseteq Z$ that is the vertex set of a component of $G[Z]$, we define $(t,C)$ to be a core; this defines the set of all cores with birthday $t$.
	For each core $(t,C)$, its {\em spread} $S(t,C)$
	is the set of all $t'\in V(T)$ such that 
	\begin{itemize}
		\item $t\le_T t'$;
		\item $C\cap B_{t'}\ne \emptyset$; and
		\item $t'\in S(s,C')$ for every core $(s,C')$ such that $s<_T t$ and $t\in S(s,C')$.
	\end{itemize}
	This completes the inductive definition of the set of all cores. 
	It is easy to check that we could replace the third bullet above with the equivalent condition:
	\begin{itemize}
		\item $C'\cap B_{t'}\ne \emptyset$ for every core $(s,C')$ such that $s<_T t$ and $t\in S(s,C')$.
	\end{itemize}
	We remark that if $(t,C),(t,C')$ are distinct cores with the same birthday then $C\cap C'=\emptyset$, but if $(t,C),(t',C')$ are cores with 
	$t\ne t'$, it is possible that $C\cap C'\ne \emptyset$, or even $C=C'$. Here are a few observations about spreads that may aid the reader's intuition:
	\begin{itemize}
		\item For each path $P$ of $T$ with one end the root $r$, and for each core $(t,C)$ with $t\in V(P)$, the set of vertices
		of $P$ that belong to the spread of $(t,C)$ depends only on which vertices $t\in V(P)$ satisfy $C\cap B_t\ne \emptyset$, and on
		which other cores have birthday in $P$ and the intersection
		of their spreads with $V(P)$. So to understand cores and spreads, it suffices to understand them when $T$ is a path.
		\item If $(t_1,C_1)$ and $(t_2,C_2)$ are cores with different birthdays, then either their spreads are disjoint, or 
		one includes the other. (We cannot say much about the spreads of cores with the same birthday.)
		\item For each core $(t,C)$, its spread is a subset of 
		$$\{t'\in V(T):t\le_T t' \text{ and }C\cap B_{t'}\ne \emptyset\}$$ 
		(which is itelf the vertex 
		set of a subtree). Indeed, the spread of $(t,C)$ is the intersection of this set and all previously-defined spreads 
		that contain $t$.
	\end{itemize}
	
	Two subsets $X,Y\subseteq V(G)$ are {\em anticomplete} if they are disjoint and there are no edges of $G$ between them.
	We need, first:
	\\
	\\
	(1) {\em If $(t_1,C_1),(t_2,C_2)$ are distinct cores and their spreads intersect, then $C_1,C_2$ are anticomplete.}
	\\
	\\
	We may assume that $t_1\ne t_2$. Since the spreads of 
	$(t_1,C_1),(t_2,C_2)$ intersect, $t_1,t_2$ have a common descendant $t_0$ say, 
	so one of $t_1,t_2$ is a strict ancestor of the other. Hence
	we may assume that $t_1<_T t_2$, 
	and therefore $t_2\in S(t_1,C_1)$ since the spreads intersect. Since
	$(t_2,C_2)$ is a core, it follows that for each $v\in C_2$, $v\notin C_1$ and $v$ has no neighbour in $C_1$.
	Consequently, $C_1,C_2$
	are anticomplete. This proves (1).
	\\
	\\
	(2) {\em For each $t\in V(T)$, there are at most $k$ cores $(s,C)$ such that $t\in S(s,C)$.}
	\\
	\\
	Let $\{(s_1,C_1)\LL (s_n,C_n)\}$ be the set of all cores whose spread contains $t$. For $1\le i\le n$, the set $C_i\cap B_t$ is nonempty; choose $v_i\in C_i\cap B_t$.
	By (1), $v_1\LL v_n$ are distinct and pairwise nonadjacent, and so $n\le \alpha(B_t)\le k$.
	This proves (2).
	
	\bigskip
	For each $v\in V(G)$, there exists $t\in V(T)$ with $v\in B_{t}$, and the set of such vertices $t$ induces a subtree of $T$.
	In particular, there is a unique $t\in V(T)$ of minimum height with $v\in B_{t}$, and we call $t$ the {\em source} of $v$.
	Let $t$ be the source of $v$. There might or might not exist $C\subseteq B_t$ with $v\in C$ such that $(t,C)$ is a core.
	If there exists such $C$ we say $v$ is {\em central}.  If there exists a core $(t',C')$
	such that $t'<_T t$ and $t\in S(t',C')$ and  $v$ has a neighbour in $C'$, we say $v$ is {\em peripheral}.
	(Note that $v$ cannot belong to $C'$, from the definition of $t$.)
	\\
	\\
	(3) {\em Every vertex $v\in V(G)$ is central or peripheral, and not both.}
	\\
	\\
	Let $t$ be the source of $v$. The first statement is clear from the definition of the set of cores with birthday $t$. For the ``not both'' part, suppose that $v$ is central and peripheral;
	choose $C\subseteq B_t$ with $v\in C$ such that $(t,C)$ is a core, and choose a core $(t',C')$ such that 
	$t'<_T t$ and $t\in S(t',C')$ and  $v$  has a neighbour in $C'$.  Since 
	$t\in S(t,C)\cap S(t', C')$, and $v\in C$ has a neighbour in $C'$, this contradicts (1), and so proves (3).
	
	\bigskip
	
	For each $v\in V(G)$, we define a set $\phi(v)$ of cores as follows. Let $t_1\in V(T)$ be the source of $v$.
	If $v$ is central, let $\phi(v)=\{(t_1,C_1)\}$ where $(t_1,C_1)$ is the core with $v\in C_1$. Now assume $v$ is peripheral.
	Hence there is a strict ancestor $t_0$ of $t_1$ and a core $(t_0,C_0)$ such that $t_1\in S(t_0,C_0)$, and $v$ has a neighbour in $C_0$.
	Choose such $t_0$ of minimum height, and let $\phi(v)$ be the set of all cores $(t_0,C_0)$
	such that $t_1\in S(t_0,C_0)$ and $v$ has a neighbour in $C_0$. 
	\\
	\\
	(4) {\em Let $t\in V(T)$ and $v\in B_t$.  Then there is a core $(t_1,C_1)$ with $t\in S(t_1,C_1)$ such that either $v\in C_1$, or $(t_1,C_1)\in \phi(v)$.}
	\\
	\\
	We proceed by induction on the height of $t$. We may assume that there is no core $(t,C)$ with $v\in C$, 
	and so $v$ was disqualified when the cores with birthday $t$ were defined. Consequently, there is a core $(t_1,C_1)$ such 
	that $t_1<_T t$, and $t\in S(t_1,C_1)$, and $v$ belongs to or
	has a neighbour in $C_1$. Choose $(t_1,C_1)$ such that $t_1$ has minimum height. 
	
	Let $s$ be the parent of $t$.
	If $v\notin B_s$, then $t$ is the source of $v$, and $(t_1,C_1)\in \phi(v)$, and the claim holds. Thus, we may assume that $v\in B_s$.
	From the inductive hypothesis, there is a core $(t_2,C_2)$ with $s\in S(t_2,C_2)$ such that either $v\in C_2$, or $(t_2,C_2)\in \phi(v)$.
	We may assume that $t\notin S(t_2,C_2)$, and so $(t_1,C_1)\ne (t_2,C_2)$.
	
	But $s$ belongs to both $S(t_1,C_1),S(t_2,C_2)$. Since $(t_1,C_1)\ne (t_2,C_2)$, it follows from (1) that $C_1,C_2$ are anticomplete.
	Since $v$ belongs to or has a neighbour in each of $C_1,C_2$, it follows that $v\notin C_1, C_2$, and so $(t_2,C_2)\in \phi(v)$. 
	Since
	$t\in S(t_1,C_1)$, the definition of $\phi(v)$ implies that $t_1\not <_T t_2$, and so $t_2\le_T t_1$.
	In particular, $t_1\in S(t_2,C_2)$; and since $t\in S(t_1,C_1)$ and $t\notin S(t_2,C_2)$, the definition of $S(t_1,C_1)$ implies that
	$t_2\not <_T t_1$. 
	Thus $t_1=t_2$, and so $(t_1,C_1)\in \phi(v)$ and the claim holds. This proves (4).
	\\
	\\
	(5) {\em Let $P$ be a path of $T$ with one end $r$, and let $v\in V(G)$. Let $t_0$ be the common birthday of the members of $\phi(v)$.
		Let $\mathcal{C}(P,v)$ be the set of cores $(t,C)$
		such that
		$t\in V(P)$ and $v\in C$. Let the members of $\mathcal{C}(P,v)$ with birthday different from $t_0$ be 
		$(t_1,C_1)\LL (t_n,C_n)$,
		numbered such that $t_0,t_1\LL t_n$ have strictly increasing height. Then:
		\begin{itemize}
			\item  for $1\le i\le n$, $t_i\notin S(t_{h}, C_{h})$ for  $1\le h<i$, and $t_i\notin  S(t_{0}, C_{0})$ for each $(t_0,C_0)\in \phi(v)$;
			\item for $1\le i\le n$, let $s_i$ be the parent of $t_i$: then $s_i\in S(t_{i-1}, C_{i-1})$ for $2\le i\le n$, and 
			there exists $(t_0,C_0)\in \phi(v)$ such that $s_1\in S(t_{0}, C_{0})$; and
			\item $n\le k$.
		\end{itemize}
	}
	\noindent
	The first bullet holds by (1), since $v\in C_i$ and either $v\in C_{h}$, or $v$ has a neighbour in $C_{0}$.
	
	For the second bullet,
	let $1\le i\le n$. If $v\notin B_{s_i}$, then $t_i$ is the source of $v$, and so $v$ is central, contradicting that $t_i\ne t_0$.
	So $v\in B_{s_i}$.
	If $s_i\in S(t_{0}, C_{0})$ for some $(t_0,C_0)\in \phi(v)$, then  $i=1$ (because otherwise $t_{i-1}\in S(t_{0}, C_{0})$ contrary to the 
	first bullet), 
	and the claim is true. So we assume that there is no $(t_0,C_0)\in \phi(v)$ with $s_i\in S(t_{0}, C_{0})$.
	From (4), 
	there is a core $(t', C')$ with $s_i\in S(t', C')$ and $v\in C'$. Hence $(t',C') = (t_h,C_h)$ for some $h\in \{1\LL i-1\}$. 
	If $h<i-1$, then $t_{i-1}\in S(t_h,C_h)$, contradicting the first bullet.
	Thus $h=i-1$. This proves the second bullet.
	
	For the third bullet, we may assume that $n\ge 1$. 
	For $1\le i\le n$, let $A_i$ be the set of cores $(t,C)$ such that 
	$t<_T t_i$
	and $t_i\in S(t,C)$. We claim that $A_i\subseteq A_{i-1}$ for $2\le i\le n$. To see this, suppose that $(t,C)\in A_i$,
	and $(t,C)\notin A_{i-1}$. Thus, $t<_T t_i$, and $t \not <_T t_{i-1}$, and consequently $t_{i-1}\le_T t$.
	There are two cases, depending whether $t=t_{i-1}$ or not. 
	\begin{itemize}
		\item Suppose first that $t_{i-1} <_T t$. 
		Then $t\in S(t_{i-1},C_{i-1})$ (because $s_i\in S(t_{i-1},C_{i-1})$ and $t\ne t_i$). But $t_i\notin S(t_{i-1}, C_{i-1})$, and $t_i\in S(t,C)$,
		contrary to the definition of $S(t,C)$. 
		\item Now suppose that $t=t_{i-1}$.  Since
		$t_i\notin S(t_{i-1},C_{i-1})$, and $C_{i-1}\cap B_{t_i}\ne \emptyset$ (because it contains $v$, since $i\ge 2$), the definition of 
		$S(t_{i-1},C_{i-1})$ implies that 
		there is a core $(d,D)$ such that $d<_T t_{i-1}$, and $t_{i-1}\in S(d,D)$, and $t_i\notin S(d,D)$. 
		But this contradicts the definition of the spread of $(t,C)$, since $d<_T t$ and $t_i\in S(t,C)$ and $t=t_{i-1}\in S(d,D)$.
	\end{itemize}
	This proves our claim that $A_i\subseteq A_{i-1}$ for $2\le i\le n$. 
	
	Next we claim that $A_i$ is a proper subset of $A_{i-1}$ for $2\le i\le n$.
	Since
	$C_{i-1}\cap B_{t_i}\ne \emptyset$ and yet
	$t_i\notin S(t_{i-1}, C_{i-1})$, there is a core $(d,D)$ such that $d<_T t_{i-1}$, and $t_{i-1}\in S(d,D)$,
	and $t_i\notin S(d,D)$. But then $(d,D)\in A_{i-1}\setminus A_i$, and so $A_i$ is a proper subset of $A_{i-1}$, as claimed.
	
	It follows that $|A_n|\le |A_1|-(n-1)$. By (2), there are at most $k$ cores $(s,C)$ such that $t_1\in S(s,C)$; and since $(t_1,C_1)$
	is one of these, and so are the all the members of $A_1$, it follows that $|A_1|\le k-1$. Consequently  $|A_n|\le (k-1)-(n-1)$,
	and so $n\le k$. This proves the third bullet, and so completes the proof of (5).
	
	\bigskip
	
	Next we construct a graph $J$. Its vertex set is the set of all triples $(s,t,C)$ where $(t,C)$ is a core and $s$ is in its spread.
	Consequently $t<_T s$ for all vertices $(s,t,C)$ of $J$. If $(s_1,t_1,C_1), (s_2,t_2, C_2)\in V(J)$ are distinct, 
	they are adjacent in $J$ if either:
	\begin{itemize}
		\item $s_1=s_2$ and $\dist_G(C_1,C_2)\le 3$; or
		\item $s_1, s_2$ are adjacent in $T$, and $C_1\cap C_2\cap B_{s}\ne \emptyset$, where $s\in \{s_1,s_2\}$ is the parent of the other member of this set.
	\end{itemize}
	In particular, if $(t,C)$ is a core and $s_1,s_2\in S(t,C)$ are adjacent in $T$, then $(s_1,t,C), (s_2,t,C)$ are adjacent in $J$; edges of $J$ of this kind are called {\em green} edges, and all other edges of $J$ are {\em red}.
	We will eventually show that 
	there is a $(2k+2,2k-1)$-quasi-isometry from $G$ to the graph obtained from $J$ by contracting all green edges. But first we prove some properties of $J$.
	\\
	\\
	(6) {\em $J$ has pseudo-tree-width at most $k$.}
	\\
	\\
	For each $s\in V(T)$, let $A_s$ be the set of all  $(s,t,C)\in V(J)$. Thus the sets $A_s\; (s\in V(T))$ are pairwise disjoint and have union
	$V(J)$. Let $s,t\in V(T)$ where $s$ is the parent of $t$. There may be edges of $J$ between $A_s$ and $A_t$, but we claim that there
	are at most $k$ such edges. 
	For each edge $e\in E(J)$ between $A_s,A_t$, we define $v_e\in V(G)$ as follows. Let the ends of $e$ be 
	$(s,s_1,C_1)$ and $(t,t_1, D_1)$.
	Then $C_1\cap D_1\cap B_s \ne \emptyset$; choose $v_e\in C_1\cap D_1\cap B_s$.
	We claim that $v_{e_1}, v_{e_2}$ are distinct and nonadjacent 
	for all distinct edges $e_1,e_2$ between $A_s,A_t$. To see this, let $e_i$ have ends $(s,s_i,C_i)$ and $(t,t_i, D_i)$ for $i = 1,2$. Either
	$(s_1,C_1)\ne (s_2,C_2)$ or $(t_1,D_1)\ne (t_2,D_2)$. In the first case, $C_1,C_2$ are anticomplete by (1),
	and so $v_{e_1}, v_{e_2}$ are distinct and nonadjacent. In the second case, $D_1,D_2$ are anticomplete by (1), and again it follows that $v_{e_1}, v_{e_2}$ are distinct and nonadjacent.
	Since $\alpha(B_s)\le k$, this proves that there are at most
	$k$ edges of $J$ between $A_s, A_t$.
	
	Let $f=st$ be an edge of $T$, where $s$ is the parent of $t$. From \ref{matching}, since $|A_s|, |A_t|\le k$ by (2), there is a 
	pseudo-path-decomposition $(B_1^f\LL B^f_{n(f)})$ of $J[A_s\cup A_t]$ with width at most $k$ and with $A_s\subseteq B^f_1$
	and $A_t\subseteq B^f_{n(f)}$. This defines $n(f)$, for each edge $f$ of $T$. Subdivide each edge $f\in E(T)$ $n(f)$ times, making a tree $T'$.
	Define $C_t=B_t$ for each $t\in V(T)$. For each $f=st\in E(T)$ where $s$ is the parent of $t$, let $s, u_1\LL u_{n(f)}, t$  be the vertices in order of the path formed
	by subdividing $f$, and define $C_{u_i} = B^f_i$ for $1\le i\le n(f)$. This defines a pseudo-tree-decomposition of $J$ with width at most $k$,
	and so proves (6).
	
	\bigskip
	
	For each $v\in V(G)$, define $\psi(v) = (t,t,C)$ where $(t,C)\in \phi(v)$ (choosing some member of $\phi(v)$ arbitrarily if there are more than one).
	\\
	\\
	(7) {\em Let $v\in V(G)$, and let $(t,C)$ be a core with $v\in C$. Then there is a path of $J$ between $\psi(v)$
		and $(t,t,C)$ with at most $k+1$ red edges.}
	\\
	\\
	Let $P$ be the path of $T$ between $r, t$, and define $t_0$ and $(t_1,C_1)\LL (t_n,C_n)$ as in (5). By the second bullet of (5),
	there exists $(t_0,C_0)\in \phi(v)$
	such that the parent of $t_1$ belongs to $S(t_0,C_0)$. Again, by the second bullet of (5),
	for $0\le i<n$, there is a path of $J$ from $(t_{i-1},t_{i-1}, C_{i-1})$ to $(t_i,t_i,C_i)$ in which all edges are green except 
	the last; and since $n\le k$ (again by (5)), and $(t,C)=(t_n,C_n)$, there is a path of $J$ from $(t_{0},t_{0}, C_{0})$ to $(t,t,C)$
	with at most $k$ red edges.  But $\psi(v)$ equals or is adjacent (by a red edge) to $(t_{0},t_{0}, C_{0})$. This proves (7).
	\\
	\\
	(8) {\em Let $v_1,v_2\in V(G)$ be adjacent. Then there is a path of $J$ between $\psi(v_1), \psi(v_2)$ using at most $k+2$ red 
		edges.}
	\\
	\\
	Let $\psi(v_i) = (t_i,t_i,C_i)$ for $i = 1,2$. Since
	$v_i$ belongs to or has a neighbour in $C_i$, for $i = 1,2$, and $v_1v_2\in E(G)$, it follows that $\dist_G(C_1,C_2)\le 3$. 
	
	A {\em green path} of $J$ means a path of $J$ containing only green edges. 
	Suppose that $t_2\in S(t_1,C_1)$, and consequently there is a green path of $J$ between $(t_1,t_1,C_1)$ and $(t_2,t_1,C_1)$.
	Since
	there is a (red) edge of $J$ between $(t_2,t_1,C_1)$ and $(t_2,t_2,C_2)$ (from the definition of $J$, since $\dist_G(C_1,C_2)\le 3$),
	the claim is true. Thus we may assume that  $t_2\notin S(t_1,C_1)$, and similarly $t_1\notin S(t_2,C_2)$. 
	
	Let $t_i'$ be the source of $v_i$ for $i = 1,2$.
	There exists $s\in V(T)$ with $v_1v_2\in B_s$, since
	$v_1v_2$ is an edge; and by choosing $s$ of minimum height we may assume that $s$ is the source of one of $v_1,v_2$, say $v_2$,
	and so $s= t_2'$, and $t_1'\le_T t_2'$, and $t_1\le_T t_1'$. Since  $t_2'\in S(t_2, C_2)$ and 
	$t_1\notin S(t_2,C_2)$, it follows that $t_1<_T t_2$; and so $t_1'<_T t_2$, since $t_1'\in S(t_1,C_1)$ and $t_2\notin S(t_1,C_1)$.
	We therefore have $t_1\le_T t_1'<_T t_2\le_T t_2'$.
	
	Since $t_2$ is in the path of $T$ between $t_1', t_2'$, and $v_1\in B_{t_1'}\cap B_{t_2'}$, it follows that $v_1\in B_{t_2}$.
	By (4), there is a core $(d,D)$ such that $t_2\in S(d, D)$ and either $v_1\in D$, or $(d,D)\in \phi(v_1)$. If $v_1\in D$,
	then $\psi(v_1)$ is joined to $(d,d,D)$ by a path of $J$ with at most $k+1$ red edges, by (7); and if $(d,D)\in \phi(v_1)$,
	then $\psi(v_1)$ is equal or adjacent (by a red edge) to $(d,d,D)$. In either case, $(d,d,D)$ is joined to $(t_2,d,D)$
	by a green path; and $(t_2,d,D)$ is adjacent to $(t_2,t_2,C_2)$ via a red edge, since $\dist_G(C_2,D)\le 2$ 
	(because $v_2$ belongs to or has a neighbour in both). This proves (8).
	\\
	\\
	(9) {\em For each core $(t,C)$, $G[C]$ has diameter at most $2k-1$.}
	\\
	\\
	By hypothesis, $\alpha(B_t)\le k$, and since $C\subseteq B_t$ it follows that 
	$G[C]$ has no induced path with $2k+1$ vertices. Since $G[C]$ is connected,
	it has diameter at most $2k-1$. This proves (9).
	\\
	\\
	(10) {\em If $(s_1,t_1,C_1)$ and $(s_2,t_2,C_2)$ are joined by a green path of $J$, and $v_1\in C_1$ and $v_2\in C_2$, 
		then $\dist_G(v_1,v_2)\le 2k-1$.}
	\\
	\\
	Any two vertices of $J$ joined by a green edge have the same second and third coordinates, and so $t_1=t_2$ and $C_1=C_2$.
	Consequently $v_1,v_2\in C_1$, and the result follows from (9). This proves (10).
	\\
	\\
	(11) {\em Let $v_1,v_2\in V(G)$, and suppose $P$ is a path of $J$ between $\psi(v_1), \psi(v_2)$ containing at most $n$ red 
		edges. Then $\dist_G(v_1,v_2)\le (2k+2)n+2k-1$.}
	\\
	\\
	If $n=0$ the result follows from (10), so we assume that $n\ge 1$.
	Let $P$ have ends $b_0$ and $a_{n+1}$, and let the red edges of $P$ be $a_1b_1,a_2b_2\LL a_nb_n$ in order, numbered such that
	there there is a green subpath of $P$ between $b_i,a_{i+1}$ for $0\le i\le n$. For $1\le i\le n$, define $\alpha_i, \beta_i\in V(G)$
	as follows: let $a_i = (s,t,C)$ and $b_i=(s',t',C')$ say; choose $\alpha_i\in C$ and $\beta_i \in C'$ with distance at most three 
	in $G$. (This is possible from the definition 
	of red edges.) Let $\beta_0=v_1$ and $\alpha_{n+1} = v_2$. Thus $\dist_G(\alpha_i,\beta_i)\le 3$ for $1\le i\le n$; 
	and $\dist_G(\beta_i,\alpha_{i+1})\le 2k-1$ by (10). Consequently $\dist_G(v_1,v_2)\le (2k+2)n+ 2k-1$.
	\\
	\\
	(12) {\em For each $j\in J$, there exists $v\in V(G)$ such that there is a path of $J$ between $j$ and $\psi(v)$ using at most $k+1$ red edges.}
	\\
	\\
	Let $j=(s,t,C)$, and choose $v\in C\cap B_s$. There is a green path between $j$ and $(t,t,C)$; and by (7), since $v\in C\subseteq B_t$,
	there is a path between $(t,t,C)$ and $\psi(v)$ containing at most $k+1$ red edges. This proves (12).
	
	\bigskip
	
	Let $H$ be obtained from $J$ by contracting all green edges. Thus each vertex of $H$ is formed by identifying
	all the vertices $(s,t,C)$ for a fixed core $(t,C)$, and so we can identify $V(H)$ with the set of all cores in the natural way.
	From (6), and since contraction does not increase pseudo-tree-width,  $H$ has pseudo-tree-width at most $k$, and from (8), (11), (12), 
	the function $\psi$ is a $(2k+2,2k-1)$-quasi-isometry 
	from $G$ to $H$. This proves \ref{quasitw2} and hence (with \ref{qtwconverse}) proves \ref{pseudotw}.~\bbox
	
	%%%%%%%%%%%%%%%%%%%%%%%%%%%%%%%%%%%%%%%%%%%%%%%%%%%%%%%%%%%%%%%%%%%%%%%%%%%%%%%%%%%%%%
	\section{Coarse line-width}
	
	In this section and the next we will prove a version of \ref{quasitw} for line-width, but first let us give the definitions more carefully.
	A {\em line} is a set $L$ that is linearly ordered by a relation $\le_L$. A {\em line-decomposition} of a graph $G$
	is a pair
	$(L,(B_t:t\in L))$ where $L$ is a line, and each $B_t$ is a subset of $V(G)$, satisfying:
	\begin{itemize}
		\item $G=\bigcup_{t\in L} G[B_t]$, and
		\item $B_t\cap B_{t''}\subseteq B_{t'}$ for all $t,t',t''\in L$ with $t\le_L t'\le_L t''$.
	\end{itemize}
	We define the {\em width} of such a decomposition to be the maximum of $|B_t|-1$ over all $t\in L$, if this maximum exists, and $\infty$ otherwise; and the {\em line-width} of $G$ is the minimum width of a line-decomposition.
	(This definition was used in~\cite{AS2, AS3, subtrees}, to increase the generality of some results about graphs with bounded path-width.
	It is also used in~\cite{diestel, thomas2}, where what we call ``line-decompositions'' are called ``linear decompositions''.)
	We will prove the following analogue of \ref{quasitw}:
	
	\begin{thm} \label{quasilw}
		For all $k,r$, there exist $L,C\ge 1$ such that if $G$ admits a line-decomposition $(L,(B_t:t\in L))$ such that for each $t\in L$,
		$B_t$ is the union of at most $k$ sets each with diameter at most $r$ in $G$, then
		$G$ admits an $(L,C)$-quasi-isometry to a graph with line-width at most $k$.
	\end{thm}
	This was first proved by Hickingbotham~\cite{hick}, with weaker constants. (Actually, he proved it at our request, because at that time 
	we needed the result for an application, and did not know how to use the approach of \ref{quasitw} to handle line-width.)
	We will deduce \ref{quasilw} as a consequence of the next result, but first we need some definitions. 
	If $L$ is a line, an {\em interval} of $L$ is a nonempty subset $I\subseteq L$ such that if $r,s,t\in L$ with $r\le_L s\le_L t$, and $r,t\in I$,
	then $s\in I$.
	An interval $I$ is {\em initial} if $I\ne L$ and
	for all $s,t\in L$ with $s\le_L t$, if $t\in I$ then $s\in I$. We say $I$ is a {\em final} interval if $L\setminus I$ is an initial interval.
	(If $I$ is an initial interval then $I, L\setminus I$ are
	both nonempty, and so $L\setminus I$ is an initial interval of the line obtained from $L$ by reversing its order.)
	If $I,J$ are intervals of $L$, we say $J$ {\em abuts} $I$  and $\{I,J\}$ is an {\em abutment} if $I\cap J=\emptyset$ and $I\cup J$ 
	is an interval of $L$.
	
	If $(L,(B_t:t\in L))$ is a line-decomposition of $G$, and $I$ is an interval, we define
	$$B(I) = \bigcup_{t\in I}B_t.$$
	If $I$ is an initial interval, then $B(I)\cap  B(L\setminus I)$
	is called the {\em $I$-split}; every path between $B(I),  B(L\setminus I)$ has a vertex in this set.
	We say $X$ is a {\em split} if $X$ is the $I$-split for some initial interval $I$.
	
	Let $(L,(B_t:t\in L))$ be a line-decomposition of a graph $G$.  A vertex $v\in V(G)$ is a {\em left-limit vertex} of 
	$(L,(B_t:t\in L))$ if for all $s,t\in L$ with $s\le_L t$, either $v\in B_s$ or $v\notin B_t$; and a {\em right-limit vertex} if for 
	all $s,t\in L$ with $s\le_L t$, either $v\notin B_s$ or $v\in B_t$. We say $v$ is a {\em limit vertex} if it is either a left- or right-limit vertex.
	
	Given a line-decomposition $(L,(B_t:t\in L))$ of a graph $G$ such that each bag is the union of $k$ sets with bounded diameter, our 
	goal is to construct a graph $H$ with line-width at most $k$ and a quasi-isometry $\phi$ from $G$ to $H$.   To do this, we will make use of the line-decomposition of $G$ as a
	``scaffold'' to construct $H$, building $H$ with a line-decomposition $(L,(C_t:t\in L))$, and a mapping $\phi$ that is suitably compatible with the line-decomposition of $G$.
	
	How do we ensure this?
	For each
	vertex $v\in V(G)$, let $I_v=\{t:v\in B_t\}$, and for each $w\in V(H)$ let $J_w=\{t:w\in C_t\}$. If we could arrange that $I_v\subseteq J_{\phi(v)}$ for every
	vertex of $G$, then for every edge $uv$ of $G$ we would know that some bag in the line-decomposition of $H$ contains both $\phi(u),\phi(v)$, and so we could make
	$\phi(u), \phi(v)$ adjacent in $H$ if we wanted.  For $\phi$ to be a quasi-isometry, we also need that the set of vertices mapping to any given $w\in V(H)$ has
	small diameter in $G$; and (to satisfy the last condition in the definition of a quasi-isometry) it will be convenient to demand that $\phi$ is surjective.
	However, it is not straightforward to construct such an $H$; and we will also want to carry out an induction that divides up $L$ into a sequence of intervals
	and then glues together  decompositions of these intervals.  So we adjust the conditions slightly.  First, we will map each vertex $v$ of $G$ to either one or two vertices of $H$; and we demand that $I_v$ is contained in the union of the ione or two corresponding intervals $J_w$.  Second, we will be more restrictive when
	handling limit vertices: these are only allowed to map to a single vertex of $H$ (which means that they will not end up mapping to many intervals at the gluing stage).  We capture these requirements with the notion of a striation.

	Let $(L,(B_t:t\in L))$ be a line-decomposition of a graph $G$. 
	For each $v\in V(G)$ let $I_v=\{t\in L:v\in B_t\}$.
	If $N$ is a set, $N^{(2)}$ denotes the set of one- and two-element subsets of $N$. 
	A {\em striation} for $(L,(B_t:t\in L))$ is a 
	pair $(\sigma, (J_n:n\in N))$, where 
	$(J_n:n\in N)$ is a family of intervals of $L$, and $\sigma$ is a map from $V(G)$ to $N^{(2)}$, satisfying:
	\begin{itemize}
		\item for each $v\in V(G)$, $I_v\subseteq \bigcup_{n\in \sigma(v)} J_n$; and moreover, if $m,n\in \sigma(v)$ are distinct, then $v$ is not a limit vertex of $(L,(B_t:t\in L))$, and
		$\{J_m, J_n\}$ is an abutment, and $I_v\cap J_m, I_v\cap J_n\ne \emptyset$;
		\item for each $n\in N$ and each $t\in J_n$, there exists $v\in V(G)$ such that $v\in B_t$ and $n\in \sigma(v)$ (in other words,
		$J_n\subseteq \bigcup_{\sigma(v)\ni n}I_v$); and
		\item for each $n\in N$, the set $\{v\in V(G):n\in \sigma(v)\}$ has diameter at most three.
	\end{itemize} 
	If $B\subseteq V(G)$ is a union of finitely many cliques of $G$, then we write $\kappa(B)$ for the minimum number of cliques in such a union (that is, the chromatic number of the complement of $G[B]$).  We say that
	a collection of
	$\kappa(B)$ cliques with union $B$ is a {\em clique min-cover} of $B$.
	We will prove:
	
	%%%%%%%%%%%%%%%%%%%%%%%%%%%%%%%%%%%
	\begin{thm}\label{getstriation}
		Let $(L,(B_t:t\in L))$ be a line-decomposition of a graph $G$, and let $k\ge 1$ such that $\kappa(B_t)\le k$ for each $t\in L$.
		Then there is a striation $(\sigma, (J_n:n\in N))$ for 
		$(L,(B_t:t\in L))$, such that $|\{n\in N: t\in J_n\}|\le \kappa(B_t)$ for each $t\in L$.
	\end{thm}
	\Proof We proceed by induction on $k$. Thus, either $k=1$ or we may assume the result holds for $k-1$. 
	\\
	\\
	(1) {\em We may assume that $G$ is connected.}
	\\
	\\
	Let $\mac C$ be the set of components of $G$. 
	For each $C\in \mac C$, $(L,(B_t\cap V(C):t\in L))$ is a line-decomposition of $C$ of finite width. Suppose that for 
	each $C$, there is a striation $(\sigma^C, (J_n:n\in N^C))$ for                   
	$(L,(B_t\cap V(C):t\in L))$, such that for each $t\in T$, the number of $n\in N^C$ such that 
	$|\{n\in N^C: t\in J_n^C\}|\le \kappa(B_t\cap V(C))$ for each 
	$t\in J_n^C$. We may assume that the sets $N^C\;(C\in \mac C)$ are pairwise 
	disjoint. Let $N$ be their union.
	For each $v\in V(G)$, define $\sigma(v) = \sigma^C(v)$ where $v\in V(C)$. Then $(\sigma, (J_n:n\in N))$ is a 
	striation for $(L,(B_t:t\in L))$. Moreover, let $t\in L$. Then for each $t\in L$,
	$$|\{n\in N: t\in J_n\}|=\sum_{C\in \mac C}|\{n\in N^C: t\in J_n^C\}|\le \sum_{C\in \mac C}\kappa(B_t\cap V(C))= \kappa(B_t),$$
	since each clique of a clique min-cover of $B$ is a subset of $B_t\cap V(C)$ for a unique $C\in \mac C$.
	Consequently, if the theorem holds for each component of $G$ then it holds for $G$.
	This proves (1).
	
	\bigskip
	
	Since $G$ is connected, it follows that if $t\in L$ and $B_t=\emptyset$, then either $B_s=\emptyset$ for all $s\in L$ with $s\le_L t$,
	or $B_s=\emptyset$ for all $s\in L$ with $t\le_L s$. Consequently the set $L'$ of $t\in L$ such that $B_t\ne \emptyset$ is an interval
	of $L$, and $(L',(B_t:t\in L'))$ is a line-decomposition of $G$ ; and if the result holds for $(L',(B_t:t\in L'))$ then it holds for
	$(L,(B_t:t\in L))$. We may therefore assume that $B_t\ne \emptyset$ for each $t\in L$. 
	\\
	\\
	(2) {\em Let $L'$ be an interval of $L$, and let $U=\bigcup_{t\in L'}B_t$.  Suppose that for each $t\in L'$, some vertex in $B_t$
		is a left-limit vertex of $(L',(B_t:t\in L'))$. Then there is a striation $(\sigma, (J_n:n\in N))$ for
		$(L',(B_t:t\in L'))$, such that $|\{n\in N: t\in J_n\}|\le \kappa(B_t)$
		for each $t\in L'$.}
	\\
	\\
	Let $S$ be the set of left-limit vertices of $(L',(B_t:t\in L'))$. For every finite subset $P$ of $S$, there exists $t\in L'$
	with $P\subseteq B_t$, and consequently $\kappa(P)\le k$ for each finite subset $P$ of $S$. By compactness (or more exactly, by a theorem of De Bruijn and Erd\H{o}s~\cite{debruijn} applied to the complement of $G[S]$), $\kappa(S)\le k$. 
	Fix some clique min-cover $\mac S$ of $S$, and choose $X\in \mac S$  with
	$\{t\in L':X\cap B_t\ne \emptyset\}$
	maximal. Thus $X$ is a clique. 
	
	We claim that $X\cap B_t\ne \emptyset$ for each $t\in L'$.
	For each $v\in U$, let $I'_v$ be the set of all $t\in L'$ with $v\in B_t$. For each $s\in S$, since $s$
	is a left-limit vertex of $(L',(B_t:t\in L'))$, the set $I_s'$
	is either null, or equal to $L'$, or an initial interval of $L'$.
	Suppose that $X\cap B_t=\emptyset$ for some $t\in L$. Then for each $s\in X$, $I'_s$ is either an initial interval of $L'$ not containing $t$, or null. 
	Let $I'_X=\bigcup_{s\in X}I'_s$. It follows that $I'_X$ is also either an initial interval of $L'$ not containing $t$, or null. Choose $s\in S$ with $s\in B_t$. Thus $I'_s$ is either an initial interval of $L'$ containing $t$, or equals $L'$,
	and in either case properly contains $I'_X$. But $s$ belongs to a member of $\mac S$, contrary to the choice of $X$. 
	This proves that $X\cap B_t\ne \emptyset$ for each $t\in L'$.
	
	Let $Y$ be the set of vertices in $U\setminus X$ with no neighbour in $X$. 
	Then $(L', (B_t\cap Y:t\in L'))$
	is a line-decomposition of $G[Y]$. Moreover, for each $t\in L'$, since $X\cap B_t\ne \emptyset$, 
	it follows that $\kappa(Y\cap B_{t})\le k-1$. If $k=1$
	then $Y=\emptyset$, and we may define $N=\{1\}$ say, and
	$J_1 = L'$,
	and define
	$\sigma(v) = 1$ for each $v\in U$, to satisfy (2). Thus we may assume that $k\ge 2$. From the inductive hypothesis,
	there is a striation $(\sigma', (J_n':n\in N'))$ for
	$(L',(B_t\cap Y:t\in L'))$,  such that $|\{n\in N': t\in J_n'\}|\le \kappa(B_t\cap Y)$ for each $t\in L'$. 
	We may assume that $1\notin N'$; let $N=N'\cup \{1\}$. For each $n\in N$, define $J_n = J'_n$
	if $n\in N'$, and $J_1 = L'$. For each $v\in U$, define $\sigma(v) = \sigma'(v)$
	if $v\in Y$, and $\sigma(v) = \{1\}$ if $v\notin Y$. We claim that
	$(\sigma, (J_n:n\in N))$ satisfies (2). To check this, we must show that:
	\begin{itemize}
		\item for each $v\in U$, $I_v'\subseteq \bigcup_{n\in \sigma(v)} J_n$; and moreover, if $m,n\in \sigma(v)$ are distinct, then 
		$v$ is not a limit vertex of $(L',(B_t:t\in L'))$, and
		$\{J_m, J_n\}$ is an abutment, and $I_v'\cap J_m, I_v'\cap J_n\ne \emptyset$;
		\item for each $n\in N$ and each $t\in J_n$, there exists $v\in U$ such that $v\in B_t$ and $n\in \sigma(v)$ (in other words,
		$J_n\subseteq \bigcup_{\sigma(v)\ni n}I_v'$);
		\item for each $n\in N$, the set $\{v\in U:n\in \sigma(v)\}$ has diameter at most three; and
		\item for each $t\in L'$, $|\{n\in N: t\in J_n\}|\le \kappa(B_t\cap U)$.
	\end{itemize}
	For the first, let $v\in U$. If $v\notin Y$, then $\bigcup_{n\in \sigma(v)} J_n=J_1=L'$, 
	and therefore includes $I_v'$.
	If $v\in Y$, 
	then $\bigcup_{n\in \sigma(v)} J_n= \bigcup_{n\in \sigma'(v)} J_n'$,
	and the latter contains $I_v'$ since
	$(\sigma', (J_n':n\in N'))$ is a striation. Moreover, if $m,n\in \sigma(v)$ are distinct, then $m,n\in \sigma'(v)$, and so 
	$v$ is not a limit vertex of $(L', (B_t\cap Y:t\in L'))$, and $\{J_m', J_n'\}$ is an abutment, and $I_v'\cap J_m', I_v'\cap J_n'\ne \emptyset$; but then 
	$v$ is not a limit vertex of $(L', (B_t:t\in L'))$, and $\{J_m, J_n\}$ is an abutment.
	
	To see the second, if $n = 1$ then $J_n=L'$, and $n\in \sigma'(v)$ if and only if $v\in U\setminus Y$;
	so we must show that $L'$ is
	the union of the sets $I_v'$ over all $v\in U\setminus Y$. But for each $t\in L'$, $X\cap B_t\ne\emptyset$,
	and so there exists $v\in U\setminus Y$ (indeed, $v\in X$) such that $t\in I_v'$, as required. If $n\ne 1$,
	then $n\in N'$, and $J_n = J'_n$, and $v\in U$ satisfies $n\in \sigma(v)$ if and only if $v\in Y$ and $n\in \sigma'(v)$.
	So we need to show that $J'_n$ is the union of the sets $I_v'$ where $v\in Y$ satisfies $n\in \sigma'(v)$, which is true
	since $(\sigma', (J_n':n\in N'))$ is a striation.
	
	To see the third, if $n=1$ then $\{v\in U:n\in \sigma(v)\}$ equals $U\setminus Y$, which has diameter at most three, since every vertex
	in $U\setminus Y$ belongs to or has a neighbour in the clique $X$. If $n\ne 1$, then $\{v\in U:n\in \sigma(v)\}$ equals
	$\{v\in Y:n\in \sigma'(v)\}$, which has diameter at most three since $(\sigma', (J_n':n\in N'))$ is a striation.
	
	Finally, for the fourth, let $t\in L'$.  Then for each $t\in L'$, 
	$$|\{n\in N: t\in J_n\}|\le |\{n\in N': t\in J_n'\}|+1,$$
	because $|N\setminus N'|=1$, and
	$J'_n = J_n$ for all $n\in N'$.  But 
	$$|\{n\in N': t\in J_n'\}|\le \kappa(B_t\cap Y),$$
	since $(\sigma', (J_n':n\in N'))$ is a striation. And 
	$$\kappa(B_t\cap Y)\le \kappa(B_t)-1,$$ 
	since if $\mac S$ is a clique min-cover of $B_t$ then at least one of its members contains a vertex in $X$ and so is disjoint from $Y$.
	This proves the fourth statement,
	and so proves (2).
	
	\bigskip
	Let $Z_1,Z_2$ be respectively the sets of left- and right-limit vertices of $(L,(B_t:t\in L))$.
	\\
	\\
	(3) {\em If $B_t\cap (Z_1\cup Z_2)\ne \emptyset$ for each $t\in L$, then the theorem holds.}
	\\
	\\
	If $B_t\cap Z_1\ne \emptyset$ for each $t\in L$, then the theorem holds by (2), so we may assume that $Z_2\ne \emptyset$ and similarly $Z_1\ne \emptyset$.
	Let $I_1$ be the set of $t\in L$ such that $B_t\cap Z_1\ne \emptyset$, and $I_2=L\setminus I_1$. 
	If $I_1=L$ then the theorem holds by (2), so we assume that $I_1\ne L$, and similarly $I_1\ne \emptyset$.
	Thus $I_1$ is an initial interval. Let $S$ be the $I_1$-split. 
	Since $Z_1\cap B_t\ne \emptyset$ for all $t\in I_1$, (2) (taking $L'=I_1$) implies that 
	there is a striation $(\sigma^1, (J_n:n\in N^1))$ for
	$(I_1,(B_t:t\in I_1))$, such that $|\{n\in N^1: t\in J_n\}|\le \kappa(B_t)$ for each $t\in I_1$.
	Since $B_t\cap Z_2\ne \emptyset$ for each $t\in I_2$, (2) implies that
	there is a striation $(\sigma^2, (J_n:n\in N^2))$ for
	$(I_2,(B_t:t\in I_2))$, such that $|\{n\in N^2: t\in J_n\}|\le \kappa(B_t)$
	for each $t\in I_2$.
	We may assume that $N^1\cap N^2=\emptyset$; let $N=N^1\cup N^2$. For $i = 1,2$ let $U_i=\bigcup_{t\in I_i}B_t$.
	Thus $B(I_1)\cap B(I_2) = S$, and $Z_1\subseteq B(I_1)$, and $Z_2\subseteq B(I_2)$. For each $v\in V(G)$, define $\sigma(v)$ as follows. 
	If $v\in B(I_i)\setminus S$ for some $i\in \{1,2\}$, let $\sigma(v) = \sigma^i(v)$. If $v\in S$, then $v$ is a limit vertex of 
	both $(I_1,(B_t:t\in I_1))$ and $(I_2,(B_t:t\in I_2))$; let $\sigma(v) = \sigma^1(v)\cup \sigma^2(v)$. 
	
	We claim that $(\sigma, (J_n:n\in N))$ satisfies the theorem. To see this, note that if $v\in S$, then $v$ is a right-limit vertex of 
	$(I_1,(B_t:t\in I_1))$, so $\sigma^1(v)$ is a singleton $\{m\}$ say, and $J_m$ either equals $I_1$ or is a final interval of $I_1$ (because 
	$I_v\cap I_1\subseteq J_m$). Similarly, $\sigma^2(v)$ is a singleton $\{n\}$ say, and $J_m$ either equals $I_2$ or is an initial 
	interval of $I_2$. Thus $\{J_m, J_n\}$ is an abutment. The remainder of the (several) things to check are easy and we omit them. (We will write them 
	out in detail after the proof of (5).) This proves (3).
	
	\bigskip
	
	Our next step is to partition $L$ into a (finite or infinite) sequence of intervals where we can apply (3).  We will use this to construct a striation; but this means we need to keep track of vertices that belong to more than one interval, and so appear in the splits between intervals.
	
	By (3), we may assume that some $B_t$ is disjoint from $Z_1\cup Z_2$. 
	An {\em integer interval} is a nonempty subset $K$ of the set $\mathbb{Z}$ of all integers, such that if $i_1\le i_2\le i_3$ are integers with
	$i_1,i_3\in K$ then $i_2\in K$. 
	\\
	\\
	(4) {\em There is an integer interval $K$ and for each $i\in K$, an initial interval $J_i$ of $L$, satisfying the following, where $S_i$ denotes the $J_i$-split:
		\begin{itemize}
			\item if $i,j\in K$ with $i<j$ then $J_i\subsetneq J_j$;
			\item for each $i\in K$, $S_i\cap (Z_1\cup Z_2)=\emptyset$, and if $i+1\in K$, then $S_i\cap S_{i+1}=\emptyset$;
			\item if $i\in K$ and $i+1\in K$, there is a clique $X$ included in one of $S_i, S_{i+1}$ such that $X\cap B_t\ne \emptyset$ for all $t\in J_{i+1}\setminus J_i$;
			\item if $i\in K$,  and $i+1\notin K$, then $B_t\cap (S_i\cup Z_2)\ne \emptyset$ for all $t\in L\setminus J_i$; and
			\item if $i\in K$ and $i-1\notin K$, then $B_t\cap (Z_1\cup S_{i})\ne \emptyset$ for all $t\in J_{i}$.
		\end{itemize}
	}
	\noindent
	We first define a set $K'$, which will eventually be the non-negative members of the set $K$ we need to construct.  
	$K'$ will consist either of all the non-negative integers, or of a finite interval of them including zero.
	Choose an initial interval $J_0$ such that the $J_0$-split is disjoint from $Z_1\cup Z_2$.
	(This is possible since some $B_t$ is disjoint from $Z_1\cup Z_2$.)
	Inductively, suppose we have chosen initial intervals 
	$J_0\LL J_{\ell}$ of $L$ for some integer $\ell\ge 0$, with corresponding splits $S_0\LL S_{\ell}$, such that 
	\begin{itemize}
		\item if $0\le i<j\le \ell $ then $J_i\subsetneq J_j$;
		\item for $0\le i\le \ell$, $S_i\ne \emptyset$, and $S_i\cap (Z_1\cup Z_2)=\emptyset$;
		\item for $0\le i<\ell$, $S_i\cap S_{i+1}=\emptyset$; and 
		\item if $0\le i<\ell$, there is a clique $X\subseteq S_i$ such that $X\cap B_t\ne \emptyset$ for all $t\in J_{i+1}\setminus J_i$.
	\end{itemize}
	Since every finite subset of $S_{\ell}$ is a subset of $B_t$ for some $t\in L$, it follows (again from the theorem of~\cite{debruijn} 
	applied in the complememt) that $\kappa(S_{\ell})\le k$.
	Let $\mac C$ be a clique min-cover of $S_{\ell}$.
	If $(S_{\ell}\cup Z_2)\cap B_t\ne \emptyset$ for all $t\in L\setminus J_{\ell}$,
	the induction is complete; so we
	assume that for each $C\in \mac C$, there exists  $t\in L\setminus J_{\ell}$ with $(C\cup Z_2)\cap B_t= \emptyset$. 
	For each $C\in \mac C$, let $M_C$ be the set of $t\in L$ such that either $t\in J_{\ell}$ or $C\cap B_t\ne \emptyset$. It follows that 
	each $M_C$ is an initial interval; and since $\mac C$ is finite and nonempty, we can choose $C\in \mac C$ with
	$M_C$ maximal. Define $J_{\ell+1} = M_C$, and let $S_{\ell+1}$ be the $J_{\ell+1}$-split. Since there exists 
	$t\in L\setminus J_{\ell}$ such that 
	$(S_{\ell}\cup Z_2)\cap B_t= \emptyset$, it follows that $t\notin J_{\ell+1}$, and so $S_{\ell+1}\cap (Z_1\cup Z_2) = \emptyset$.  
	Moreover, 
	$C\cap B_t\ne \emptyset$ for all $t\in J_{\ell+1}\setminus J_{\ell}$. We claim that 
	$S_{\ell}\cap S_{\ell+1}=\emptyset$. Suppose not, and let $v\in S_{\ell}\cap S_{\ell+1}$. Choose $C'\in \mac C$
	with $v\in C'$, 
	%Let $C'$ be the component of $G[S_{\ell}]$ that contains $v$. 
	and choose
	$t\in L\setminus M_C$ with $v\in B_t$. It follows that $t\in M_{C'}$, and so $t\in M_C$ from the maximality of $M_C$, contradicting
	that $t\notin M_C$. This proves that $S_{\ell}\cap S_{\ell+1}=\emptyset$, and so completes the inductive step. 
	
	Let $K'$ be the set of all the integers $i\ge 0$ such that $J_i$ is defined in the process above.
	Thus $K'$ consists either of all the non-negative integers, or of $\{0\LL \ell\}$ for some integer $\ell\ge 0$.
	In the second case, there is a clique $X$ such that $X\cap B_t\ne \emptyset$ for all $t\in L\setminus J_{\ell}$.
	Consequently, so far we have constructed an integer interval $K'$ with $0\in K'$, containing only non-negative integers,
	such that 
	\begin{itemize}
		\item if $i,j\in K'$ with $i<j$ then $J_i\subsetneq J_j$;
		\item for each $i\in K'$, $S_i\ne \emptyset$, and $S_i\cap (Z_1\cup Z_2)=\emptyset$;
		\item if $i, i+1\in K'$ then $S_i\cap S_{i+1}=\emptyset$;
		\item if $i\in K'$ and $i+1\in K'$, there is a clique $X\subseteq S_i$ such that $X\cap B_t\ne \emptyset$ for all $t\in J_{i+1}\setminus J_i$; and
		\item if $i\in K'$ and $i+1\notin K'$, then  $(S_{i}\cup Z_2)\cap B_t\ne \emptyset$ for all $t\in L\setminus J_{i}$.
	\end{itemize}
	Now we repeat the process, inductively constructing $J_{-1}, J_{-2}$ and so on, with the same argument, used with the order of $L$ reversed. This proves (4).
	
	\bigskip
	If $i\notin K$ and $i+1\in K$, let $J_i=\emptyset$, and if $i\in K$ and $i+1\notin K$ let $J_{i+1} = L$. 
	Let $K^+$ be the set of all integers $i$ such that at least one of $i,i+1$ is in $K$. For each 
	integer $i\in K^+$, let $L^i = J_{i+1}\setminus J_i$, with order inherited from $L$, and let $U^i=\bigcup_{t\in L^i}B_t$.
	\\
	\\
	(5) {\em For each $i\in K^+$, there is a striation $(\sigma^i, (J^i_n:n\in N^i))$ for
		$(L^i,(B_t:t\in L^i))$ such that $|\{n\in N^i:t\in J_n^i\}|\le \kappa(B_t)$ for each $t\in L$.}
	\\
	\\
	The third, fourth or fifth bullet of (4) (whichever applies) shows
	that for each $t\in L^i$, $B_t$ contains a limit vertex of
	$(L^i,(B_t:t\in L^i))$.
	Hence the claim follows from (3) applied to $L^i$. This proves (5).
	
	\bigskip
	
	We may assume that the sets $N^i\;(i\in K^+)$ are pairwise disjoint. Let $N$ be their union. For each $n\in N$, let $J_n=J^i_n$
	where $n\in N^i$. 
	For each $i\in K$, let $S_i$ be the $J_i$-split. For each $v\in V(G)$, we define $\sigma(v)$ as follows.
	
	If $v\notin \bigcup_{i\in K}S_i$, there is a unique value of $i$ such that $v\in B_t$ for some $t\in J_{i+1}\setminus J_i$; define $\sigma(v) = \sigma^i(v)$. If $v\in \bigcup_{i\in K}S_i$,
	there is a unique value of $i\in K$ such that $v\in S_i$ (because the sets $S_i\;(i\in K)$ are pairwise disjoint). Then $v$ 
	is a right-limit vertex of $(L^{i-1},(B_t:t\in L^{i-1}))$ and a left-limit vertex of $(L^{i},(B_t:t\in L^{i}))$. Consequently,
	$|\sigma^{i-1}(v)|=|\sigma^{i}(v)|=1$. Let $\sigma^{i-1}(v)=\{m\}$ and $\sigma^{i}(v)=\{n\}$. Since $J^i_n$ includes each 
	$t\in L^i$ with $v\in B_t$, and $v$ is a left-limit vertex of  $(L^i,(B_t:t\in L^i))$, it follows that 
	$J^i_n$ either equals $L^i$ or is an initial interval of $L^i$; and similarly $J^{i-1}_m$ either equals $L^{i-1}$ or is a final interval 
	of $L^{i-1}$. Consequently $\{J_m, J_n\}$ is
	an abutment. Let $\sigma(v) = \{m,n\}$.
	
	For each $v\in V(G)$, let $I_v = \{t\in L:v\in B_t\}$ as usual.
	We claim that $(\sigma, (J_n:n\in N))$ is a striation for
	$(L,(B_t:t\in L))$, such that $|\{n\in N:t\in J_n\}|\le \kappa(B_t)$ for each $t\in L$. 
	To show this, we must check:
	\begin{itemize}
		\item for each $v\in V(G)$, $I_v\subseteq \bigcup_{n\in \sigma(v)} J_n$; and moreover, if $m,n\in \sigma(v)$ are distinct, then $v$ is not a limit vertex of $(L,(B_t:t\in L))$, and
		$\{J_m, J_n\}$ is an abutment, and $I_v\cap J_m, I_v\cap J_n\ne \emptyset$;
		\item for each $n\in N$ and each $t\in J_n$, there exists $v\in V(G)$ such that $v\in B_t$ and $n\in \sigma(v)$; and
		\item for each $n\in N$, the set $\{v\in V(G):n\in \sigma(v)\}$ has diameter at most three; and
		\item $|\{n\in N:t\in J_n\}|\le \kappa*B_t)$ for each $t\in L$.
	\end{itemize}
	For the first, let $v\in V(G)$. Suppose first that $v\notin \bigcup_{i\in K}S_i$. Then there is a unique value of $i\in K^+$
	such that $v\in U^i$; and then $\sigma(v) = \sigma^i(v)$. Moreover, $I_v\subseteq L^i$, and so $I_v\subseteq \bigcup_{n\in \sigma(v)} J_n$
	since $(\sigma^i, (J^i_n:n\in N^i))$ is a striation. If $m,n\in \sigma(v)$ are distinct, then $v$ is not a limit vertex of 
	$(L^i,(B_t:t\in L^i))$, and therefore not a limit vertex of $(L,(B_t:t\in L))$; and $\{J_m,J_n\} = \{J^i_m,J^i_n\}$ is an abutment.
	
	Now suppose that $v\in S_i$ for some (necessarily unique) value of $i\in K$. Thus $\sigma(v) = \{m,n\}$ where
	$\sigma^i(v)=\{n\}$ and $\sigma^{i-1}(v)=\{m\}$.
	Then $\emptyset\ne I_v\cap L_i\subseteq J_n$,
	and $I_v\cap L_{i-1}\subseteq J_m$, and $I_v\subseteq L_{i-1}\cup U^i$. We have seen that $\{J_m,J_n\}$ is an abutment. Since $v\in S_i$, 
	and $S_i$ is disjoint from $Z_1\cup Z_2$, it follows that $v$ is not
	a limit vertex of $(L,(B_t:t\in L))$. This proves the first statement. 
	The second, third and fourth statements are clear. 
	
	This proves that $(\sigma, (J_n:n\in N))$  satisfies the theorem, and so proves \ref{getstriation}.~\bbox
	
	%%%%%%%%%%%%%%%%%%%%%%%%%%%%%%%%%%%%%%%%%%%%%%%%%%%%%%%%%%%%%%%%%%%%%%%%%%%%%%%%%%%%%%%%%%%%%%%%%%%%%%%%%%%%%%%%%%%%%%%%%%%%%%%%%%%
	\section{Pseudo-line-decompositions}
	If $L$ is a line, we say $s,t\in L$ are {\em consecutive in $L$} if there is no $r\in L$ with $r\ne s,t$ such that 
	$s\le_L r\le_L t$ or $t\le_L r\le_L s$.
	A {\em pseudo-line-decomposition} of a graph $G$ is a pair $(L,(B_t:t\in L))$, where $L$ is a line, and $B_t$
	is a subset of $V(G)$ for each $t\in L$,  such that:
	\begin{itemize}
		\item $V(G)=\bigcup_{t\in L}B_t$;
		\item for every edge $e=uv$ of $G$, either there exists $t\in L$ with $u,v\in B_t$, or there exist consecutive $s,t\in L$
		such that $B_s\setminus B_t = \{u\}$ and $B_t\setminus B_s = \{v\}$; and
		\item for all $t_1,t_2,t_3\in L$, if $t_1\le_L t_2\le_L t_3$ then 
		$B_{t_1}\cap B_{t_3}\subseteq B_{t_2}$.
	\end{itemize}
	The {\em width} of a pseudo-tree-decomposition $(L,(B_t:t\in L))$ is the maximum of the numbers $|B_t|$ for $t\in L$,
	or $\infty$ if there is no finite maximum;
	and the {\em pseudo-line-width} of $G$ is the minimum width of a pseudo-line-decomposition of $G$.
	
	If $\{J_1,J_2\}$ is an abutment of $L$, there is a unique initial interval $I$ of $L$ that contains one of $J_1,J_2$
	and is disjoint from the other, and we call $I$ the {\em border} of $\{J_1,J_2\}$.
	We need the following lemma, which is proved simply by applying \ref{matching} to each initial interval of $L$:
	%%%%%%%%%%%%%%%%%%%%%%%
	\begin{thm}\label{abut}
		Let $H$ be a graph, let $L$ be a line, and for each $t\in L$ let $C_t\subseteq V(H)$, such that:
		\begin{itemize}
			\item $V(H)=\bigcup_{t\in L}C_t$;
			\item for all $t_1,t_2,t_3\in L$, if $t_1\le_L t_2\le_L t_3$ then
			$C_{t_1}\cap C_{t_3}\subseteq C_{t_2}$; and
			\item  for every edge $mn\in E(H)$, either there exists $t\in L$ with $m,n\in C_t$, or 
			$\{J_m,J_n\}$
			is an abutment in $L$ (where $J_m = \{t\in L:m\in C_t\}$ and so on).
		\end{itemize}
		Suppose that $k\ge 0$ is an integer such that:
		\begin{itemize}
			\item $|C_t|\le k$ for each $t\in L$; and
			\item for every initial interval $I$ of $L$, there are at most $k$ edges $mn$ of $H$ such that 
			$\{J_m, J_n\}$
			is an abutment with border $I$.
		\end{itemize}
		Then $H$ has pseudo-line-width at most $k$.
	\end{thm}
	
	The main result of this section is:
	
	%%%%%%%%%%%%%%%%%%%%%%%%%%%%%%%%
	\begin{thm}\label{mainlinethm}
		For all $k,r$, if $G$ admits a line-decomposition $(L,(B_t:t\in L))$ such that for each $t\in L$,
		$B_t$ is the union of at most $k$ sets each with diameter at most $r$ in $G$, then
		$G$ admits a $(7,3)$-quasi-isometry to a graph with pseudo-line-width at most $k$.
		
		Conversely, for all $L,C\ge 1$,
		if $G$ admits an $(L,C)$-quasi-isometry to a graph with pseudo-line-width at most $k$,
		then $G$ admits a line-decomposition $(L,(B_t:t\in L))$ such that for each $t\in L$, $B_t$ is the union of at most
		$k$ sets each of diameter at most $2L(L+C)+C$.
	\end{thm}
	
	The second half of the theorem is proved exactly as \ref{qtwconverse}, changing ``path'' to ``line'' in the appropriate places, and 
	we omit it. Moreover, the discussion in the paragraph before the statement of \ref{quasitw2} applies equally to line-decompositions;
	so to prove \ref{mainlinethm} it suffices to prove the following:
	%%%%%%%%%%%%%%%%%%%%%%%%%%%%
	\begin{thm}\label{quasilw2}
		For all $k$, if $G$ admits a line-decomposition $(L,(B_t:t\in L))$ such that $\kappa(B_t)\le k$ for each $t\in L$,
		then
		$G$ admits a $(7,3)$-quasi-isometry to a graph with pseudo-line-width at most $k$.
	\end{thm}
	\Proof
	By \ref{getstriation}, there is a striation $(\sigma, (J_n:n\in N))$ for
	$(L,(B_t:t\in L))$, such that for each $t\in L$, $|\{n\in N:t\in J_n\}|\le \kappa(B_t)\le k$. 
	Let $H$ be a graph with vertex set $N$, defined as follows.
	\begin{itemize}
		\item Distinct $n_1,n_2\in N$ are adjacent in $H$ if there are adjacent $v_1,v_2\in V(G)$ such that 
		$n_1\in \sigma(v_1)$, $n_2\in \sigma(v_2)$, and 
		$$I_{v_1}\cap J_{n_1}\cap I_{v_2}\cap J_{n_2}\ne \emptyset.$$
		\item In addition, for every initial interval $I$ of $L$ and every component $C$ of the $I$-split, let $C'$ be the set of $v\in C$
		such that there exist $m,n\in N$ with $\sigma(v)=\{m,n\}$, 
		$J_m\subseteq I$ and $J_n\cap I=\emptyset$. For each $C$ with $C'\ne \emptyset$, 
		choose exactly one $v\in C'$ 
		(called the {\em $(I,C)$-representative.}), and make the members of $\sigma(v)$ adjacent in $H$.
	\end{itemize}
	We deduce:
	\\
	\\
	(1) {\em $H$ has pseudo-line-width at most $k$.}
	\\
	\\
	For each $t\in L$, let $C_t$ be the set of all $n\in N$ such that $t\in J_n$. It follows that $V(H)=\bigcup_{t\in L}C_t$, and
	for all $t_1,t_2,t_3\in L$, if $t_1\le_L t_2\le_L t_3$ then
	$C_{t_1}\cap C_{t_3}\subseteq C_{t_2}$.
	
	We claim that
	for every edge $n_1n_2\in E(H)$, either there exists $t\in L$ with $n_1,n_2\in C_t$, or $\{J_{n_1},J_{n_2}\}$
	is an abutment in $L$. To see this, let $n_1n_2\in E(H)$. If there are adjacent $v_1,v_2\in V(G)$ such that 
	there exists $t\in I_{v_1}\cap J_{n_1}\cap I_{v_2}\cap J_{n_2}\ne \emptyset$,  
	then $t\in J_{n_i}$ for $i = 1,2$, and so $n_1,n_2\in C_t$, as claimed. On the other hand, if there exists $v\in V(G)$ with $n_1,n_2\in \sigma(v)$, then $\sigma(v) = \{n_1,n_2\}$, and $\{J_{n_1},J_{n_2}\}$ is an abutment, as claimed.
	
	Next, we claim that 
	\begin{itemize}
		\item $|C_t|\le k$ for each $t\in L$; and
		\item for every initial interval $I$ of $L$, there are at most $k$ edges $mn$ of $H$ such that 
		$\{J_m, J_n\}$
		is an abutment with border $I$.
	\end{itemize}
	To see the first statement, let $t\in L$. From the choice of the striation, 
	the number of $n\in N$ such that $t\in J_n$ is
	at most $k$, and so $|C_t|\le k$. For the second statement, let $I$ be an initial interval of $L$ and let $S$ be the $I$-split. 
	Let $F$ be the
	set of edges $mn\in E(H)$ such that $\{J_m, J_n\}$
	is an abutment with border $I$. For each $mn\in F$, since $J_m\cap J_n=\emptyset$, $mn$ is contained in $E(H)$ because of the 
	second bullet in the definition of $H$; and so there is an initial interval $I'$ of $L$ and a component $C$ of the $I'$-split
	such that $\sigma(v)=\{m,n\}$ where $v$ is the $(I',C)$-representative. Thus $I'$ is the border of the abutment $\{J_m,J_n\}$.
	On the other hand, since $mn\in F$, it follows that $I$ is the border of $\{J_m,J_n\}$, and so $I'=I$. Since $G[S]$
	has at most $k$ components, and there is at most one $(I,C)$-representative for each component $C$, it follows that $|F|\le k$.
	From \ref{abut}, this proves (1).
	
	\bigskip
	
	For each $v\in V(G)$, let $\phi(v)$
	be some member of $\sigma(v)$. To complete the proof, we will show that $\phi$ is an $(7,3)$-quasi-isometry from $G$ to $H$. 
	\\
	\\
	(2) {\em For each $v\in V(G)$, if $m,n\in \sigma(v)$ then $\dist_H(m,n)\le 3$. Consequently,
		for each $n\in V(H)$, there exists $v\in V(G)$ such that $\dist_H(\phi(v),n)\le 3$.}
	\\
	\\
	Since $m,n\in \sigma(v)$ and we may assume that $m\ne n$, it follows that $\{J_m,J_n\}$ is an abutment; let $I$ be its border,
	and let $S$ be the $I$-split. Since $I_v\cap J_m, I_v\cap J_n\ne \emptyset$, it follows that $v\in S$; let $C$ be the component
	of $S$ containing $v$. Since $v\in V(C)$ and $\{J_m,J_n\}$ is an abutment with border $I$, it follows that there is an 
	$(I,C)$-representative $v'$ say. Thus $v,v'$ are equal or adjacent in $G$.  Let $\sigma(v') = \{m',n'\}$ where $J_{m'}\subseteq I$ 
	and $J_{n'}\cap I = \emptyset$; so $m', n'$ are adjacent in $H$. If $v'=v$ then $m'=m$ and $n'=n$ and $\dist_H(m,n)\le 1$
	as required, so we assume that $v'\ne v$. Since $v,v'\in S$, there exists $t\in I$ such that $v,v'\in B_t$. But $I_v\subseteq J_m\cup J_n$,
	and so $I_v\cap I\subseteq J_m$, and therefore $t\in I_{v}\cap J_{m}$. Similarly $t\in I_{v'}\cap J_{m'}$, and so 
	$$I_{v}\cap J_{m}\cap I_{v'}\cap J_{m'}\ne \emptyset,$$
	and therefore $m,m'$ are equal or adjacent in $H$ from the definition of $H$. Similarly $n,n'$ are equal or adjacent, and
	so $\dist_H(m,n)\le 3$. This proves the first statement of (2). 
	
	For the second statement, let $n\in V(H)$, and choose $v\in V(G)$ with $n\in \sigma(v)$ (this is possible from the definition of a striation). We may assume that
	$\phi(v)\ne n$, and so $\phi(v) = m$ where $\sigma(v) = \{m,n\}$. Consequently $\dist_H(\phi(v),n)\le 3$ by the first statement. 
	This proves (2).
	\\
	\\
	(3) {\em For all $u,v\in V(G)$, $\dist_H(\phi(u),\phi(v))\le 7\dist_G(u,v)$.}
	\\
	\\
	First, let us assume that $u,v\in V(G)$ are adjacent. 
	Since $(L, (B_t:t\in L))$ is a line-decomposition, there exists $t\in L$ with $u,v\in B_t$. Since 
	$t\in I_u\subseteq \bigcup_{m\in \sigma(u)}J_m$, there exists $m\in \sigma(u)$ with $t\in J_m$, and similarly there exists 
	$n\in \sigma(v)$ with $t\in J_{n}$.
	Consequently $\dist_H(m,n)\le 1$, since 
	$$I_{u}\cap J_{m}\cap I_{v}\cap J_{n}\ne \emptyset.$$
	But by (2), $\dist_H(\phi(u),m)\le 3$, and $\dist_H(\phi(v), n)\le 3$, and so $\dist_H(\phi(u),\phi(v))\le 7$. 
	Hence, for every choice of $u,v\in V(G)$, it follows that $\dist_H(\phi(u),\phi(v))\le 7\dist_G(u,v)$, by summing what we just proved over the edges of a shortest path of $G$ joining $u,v$. This proves (3).
	\\
	\\
	(4) {\em For all $u,v\in V(G)$, and all $m\in \sigma(u)$ and $n\in \sigma(v)$, if $m=n$ then $\dist_G(u,v)\le 3$
		and if $m,n$ are adjacent in $H$, then $\dist_G(u,v)\le 7$.}
	\\
	\\
	If $m=n$, then $\dist_G(u,v)\le 3$, since $\{w\in V(G):n\in \sigma(w)\}$ has diameter at most three from the definition of           
	a striation. Thus we assume that $m\ne n$. 
	There are two reasons why $m,n$ might be adjacent in $H$: either
	\begin{itemize}
		\item there are adjacent $v_1,v_2\in V(G)$ such that $m\in \sigma(v_1)$, $n\in \sigma(v_2)$, and there exists 
		$$t\in I_{v_1}\cap J_{m}\cap I_{v_2}\cap J_{n},$$
		or 
		\item there is an initial interval $I$ of $L$ and a component $C$ of the $I$-split, and $w\in C$ such that $w$ is the $(I,C)$-representative and
		$\sigma(w)=\{m,n\}$ where $J_m\subseteq I$ and $J_n\cap I=\emptyset$.
	\end{itemize}
	In the first case, since $m\in \sigma(u)\cap \sigma(v_1)$, it follows that $\dist_G(u,v_1)\le 3$, and similarly $\dist_G(v,v_2)\le 3$,
	and so $\dist_G(u,v)\le 7$. In the second case, since $m\in \sigma(u)\cap \sigma(w)$ it follows that $\dist_G(u,w)\le 3$, and 
	similarly $\dist_G(v,w)\le 3$, and so $\dist_G(u,v)\le 6$. This proves (4).
	\\
	\\
	(5) {\em For all $u,v\in V(G)$, $\dist_G(u,v)\le  7\dist_H(\phi(u),\phi(v))+3$.}
	\\
	\\
	Let $d=\dist_H(\phi(u),\phi(v))$. As before, we may assume that $d\ge 1$. 
	Let $P$ be a shortest path of $H$ between $\phi(u), \phi(v)$, and let its vertices be $n_0\LL n_d$ in order, where $n_0 = \phi(u)$
	and $n_d=\phi(v)$. 
	For each $i\in \{1\LL d-1\}$, choose $v_i\in V(G)$ such that $n_i\in \sigma(v_i)$; and define
	$v_0 = u$ and $v_d = v$. Hence, by (4), for $1\le i\le d$, $\dist_G(v_{i-1},v_i)\le 7$, and so $\dist_G(u,v)\le 7d$.
	This proves (5).
	
	\bigskip
	
	From (2), (3) and (5), it follows that $\phi$ is a $(7,3)$-quasi-isometry from $G$ to $H$. This proves \ref{quasilw2}.~\bbox
	
	\section*{Acknowledgement}
	Our thanks to Robert Hickingbotham, for his contributions to proving \ref{quasilw}. And we would particularly like to thank 
	a referee for detecting some significant errors.
	
	\section*{Notes}
	For the purpose of open access, the author has applied a CC BY public copyright licence to
	any author accepted manuscript arising from this submission.

\end{document}